\magnification=\magstep1
\def\nz{\hfill\break\noindent}
\def\mn{{\medskip\noindent}}
\def\bn{{\bigskip\noindent}}
\def\sn{{\smallskip\noindent}}

\catcode`@=11
\font\teneufm=eufm10
\font\seveneufm=eufm7
\font\fiveeufm=eufm5
\newfam\eufmfam
\textfont\eufmfam=\teneufm
\scriptfont\eufmfam=\seveneufm
\scriptscriptfont\eufmfam=\fiveeufm
\def\frak{%
 \let\next\relax
 \ifmmode\let\next\frak@
 \else\def\next{\errmessage{Use \string\frak\space in math mode only}}
 \fi
 \next}
\def\frak@#1{{\frak@@{#1}}}
\def\frak@@#1{\fam\eufmfam#1}
\let\goth\frak
\catcode`\@=12

\def\bbbr{{\rm I\!R}}

\def\bbbn{{\rm I\!N}}

\def\bbbt{{\rm I\!T}}

\def\bbbone{{\mathchoice {\rm 1\mskip-4mu l} {\rm 1\mskip-4mu l}
{\rm 1\mskip-4.5mu l} {\rm 1\mskip-5mu l}}}
\def\bbbc{{\mathchoice {\setbox0=\hbox{$\displaystyle\rm
C$}\hbox{\hbox
to0pt{\kern0.4\wd0\vrule height0.9\ht0\hss}\box0}}
{\setbox0=\hbox{$\textstyle\rm C$}\hbox{\hbox
to0pt{\kern0.4\wd0\vrule height0.9\ht0\hss}\box0}}
{\setbox0=\hbox{$\scriptstyle\rm C$}\hbox{\hbox
to0pt{\kern0.4\wd0\vrule height0.9\ht0\hss}\box0}}
{\setbox0=\hbox{$\scriptscriptstyle\rm C$}\hbox{\hbox
to0pt{\kern0.4\wd0\vrule height0.9\ht0\hss}\box0}}}}
\def\bbbe{{\mathchoice {\setbox0=\hbox{\smalletextfont
e}\hbox{\raise
0.1\ht0\hbox to0pt{\kern0.4\wd0\vrule width0.3pt
height0.7\ht0\hss}\box0}}
{\setbox0=\hbox{\smalletextfont e}\hbox{\raise
0.1\ht0\hbox to0pt{\kern0.4\wd0\vrule width0.3pt
height0.7\ht0\hss}\box0}}
{\setbox0=\hbox{\smallescriptfont e}\hbox{\raise
0.1\ht0\hbox to0pt{\kern0.5\wd0\vrule width0.2pt
height0.7\ht0\hss}\box0}}
{\setbox0=\hbox{\smallescriptscriptfont e}\hbox{\raise
0.1\ht0\hbox to0pt{\kern0.4\wd0\vrule width0.2pt
height0.7\ht0\hss}\box0}}}}
\def\bbbq{{\mathchoice {\setbox0=\hbox{$\displaystyle\rm
Q$}\hbox{\raise
0.15\ht0\hbox to0pt{\kern0.4\wd0\vrule height0.8\ht0\hss}\box0}}
{\setbox0=\hbox{$\textstyle\rm Q$}\hbox{\raise
0.15\ht0\hbox to0pt{\kern0.4\wd0\vrule height0.8\ht0\hss}\box0}}
{\setbox0=\hbox{$\scriptstyle\rm Q$}\hbox{\raise
0.15\ht0\hbox to0pt{\kern0.4\wd0\vrule height0.7\ht0\hss}\box0}}
{\setbox0=\hbox{$\scriptscriptstyle\rm Q$}\hbox{\raise
0.15\ht0\hbox to0pt{\kern0.4\wd0\vrule
height0.7\ht0\hss}\box0}}}}
\def\bbbt{{\mathchoice {\setbox0=\hbox{$\displaystyle\rm
T$}\hbox{\hbox to0pt{\kern0.3\wd0\vrule
height0.9\ht0\hss}\box0}}
{\setbox0=\hbox{$\textstyle\rm T$}\hbox{\hbox
to0pt{\kern0.3\wd0\vrule height0.9\ht0\hss}\box0}}
{\setbox0=\hbox{$\scriptstyle\rm T$}\hbox{\hbox
to0pt{\kern0.3\wd0\vrule height0.9\ht0\hss}\box0}}
{\setbox0=\hbox{$\scriptscriptstyle\rm T$}\hbox{\hbox
to0pt{\kern0.3\wd0\vrule height0.9\ht0\hss}\box0}}}}
\def\bbbs{{\mathchoice
{\setbox0=\hbox{$\displaystyle     \rm
S$}\hbox{\raise0.5\ht0\hbox
to0pt{\kern0.35\wd0\vrule height0.45\ht0\hss}\hbox
to0pt{\kern0.55\wd0\vrule height0.5\ht0\hss}\box0}}
{\setbox0=\hbox{$\textstyle        \rm
S$}\hbox{\raise0.5\ht0\hbox
to0pt{\kern0.35\wd0\vrule height0.45\ht0\hss}\hbox
to0pt{\kern0.55\wd0\vrule height0.5\ht0\hss}\box0}}
{\setbox0=\hbox{$\scriptstyle      \rm
S$}\hbox{\raise0.5\ht0\hbox
to0pt{\kern0.35\wd0\vrule height0.45\ht0\hss}\raise0.05\ht0\hbox
to0pt{\kern0.5\wd0\vrule height0.45\ht0\hss}\box0}}
{\setbox0=\hbox{$\scriptscriptstyle\rm
S$}\hbox{\raise0.5\ht0\hbox
to0pt{\kern0.4\wd0\vrule height0.45\ht0\hss}\raise0.05\ht0\hbox
to0pt{\kern0.55\wd0\vrule height0.45\ht0\hss}\box0}}}}

\font\tifont=cmbx8 scaled \magstep3

\def\litem{\par\noindent
               \hangindent=\parindent\ltextindent}
\def\ltextindent#1{\hbox to \hangindent{#1\hss}\ignorespaces}
\def\qed{\vrule height 1.2ex width 1.1ex depth -.1ex}
\def\={\!=\!}

\def\X1{{\tilde X}}
\def\S{{SL_q(2)}}

\def\ker{{\rm ker}\,}

\def\u{{{\cal U}_q(sl(2))}}

\def\A{{\cal A}}
\def\A1{{\tilde{A}}}
\def\B{{\cal B}}

\def\R{{\cal R}}
\def\E{{\cal E}}
\def\B{{\cal B}}
\def\U{{\cal U}}

\def\V{{\cal V}}
\def\V1{{\tilde{V}}}

\def\D{{\cal {D}}}
\def\D1{{\tilde{D}}}
\def\S1{{\cal S}}
\def\X{{\cal X}}

\def\n{{\goth{n}}}
\def\m{{\goth{m}}}
\def\i{{\goth{i}}}
\def\j{{\goth{j}}}
\def\k{{\goth{k}}}

\def\1{{\bbbone}}
\def\d{{\rm d}}

\hfill q-alg/9601020
\bigskip
\noindent

{\tifont \centerline{{\bf Left-Covariant Differential Calculi on 
${\bf SL_q(2)}$ and ${\bf SL_q(3)}$}}}
\medskip\bigskip\noindent
\centerline{Konrad Schm\"udgen, Axel Sch\"uler}
\bigskip\noindent

\centerline{Fakult\"at f\"ur Mathematik und Informatik und
Naturwiss.-Theoretisches Zentrum,}
\centerline{Universit\"at Leipzig, Augustusplatz 10,}
\centerline{ 04109 Leipzig, Federal Republic of Germany;}
\centerline{ e-mail: schmuedgen@mathematik.uni-leipzig.d400.de}
\centerline{schueler@server1.rz.uni-leipzig.de}
\bigskip\noindent
{\bf Abstract:}\nz
We study $N^2-1$ dimensional left-covariant differential
calculi on the quantum group $SL_q(N)$ for which the 
generators of the quantum Lie algebras annihilate the quantum trace.
In this way we obtain one distinguished calculus on $SL_q(2)$ (which 
corresponds
to Woronowicz' 3D-calculus on $SU_q(2)$) and two distinguished calculi 
on
$SL_q(3)$ such that the higher order calculi give the ordinary differential
calculus on $SL(2)$ and $SL(3)$, respectively, in the limit $q\rightarrow 1$.
Two new differential calculi on $SL_q(3)$ are introduced and developed
 in 
detail.\bn
{\bf MSC (1991):} 17~B~37, 46~L~87, 81~R~50.\bn
{\bf 0. Introduction}\nz 
After the seminal work [18] of S.L. Woronowicz, bicovariant
differential calculi on quantum groups (Hopf algebras) have been 
extensively
studied in the literature. There is a well developed general theory of such 
calculi. Bicovariant differential calculi on the quantum group $SL_q(N), 
N\ge3$, 
have been recently classified in [12]. The case of $SL_q(2)$ has been 
treated 
before in  [15] und in  [11]. All calculi occuring in this classification
have dimension $N^2$, i.e. their dimension does not coincide with the
dimension $N^2-1$ of the corresponding classical Lie group. On the other
hand, the first example of a non-commutative differential calculus on a 
quantum group was Woronowicz' 3D-calculus on $SU_q(2)$ [17]. This is a
 three
dimensional left-covariant calculus which is not bicovariant. The 
3D-calculus is algebraically much simpler and in many respects nearer 
to the 
classical differential calculus on $SU_q(2)$ than the four dimensional 
bicovariant calculi on  $SU_q(2)$. This motivates to look for $N^2-1$ 
dimensional left-covariant differential calculi on $SL_q(N)$.  
The purpose of this paper is to study the cases $N=2$ and $N=3$
in detail. The main aim 
of our approach is to follow the classical situation as close as possible.\nz
Let us briefly explain the basic idea of the approach given in this paper. 
As in [12], we assume 
that the differentials $\d u^i_j$ of the matrix entries $u^i_j$ generate the 
left module of 1-forms. Hence the differential $\d$ can be expressed as 
$\d x=\sum (X_{ij}\ast x)\omega (u^i_j)$ for $x\in SL_q(N)$, where 
$\omega (u^i_j):= \sum_n\kappa (u^i_n) \d u^n_j$ are the left-invariant 
Maurer-Cartan 
forms and $X_{ij}$ are linear functionals on $SL_q(N)$ such that $X_{ij} 
(1)=0$. In our approach for $N=2,3$ 
the functionals $X_{ij}$ will be chosen {}From the quantized universal 
enveloping algebra $\U_q(sl_N)$. For the functionals $X_{ij}$ with 
$i\ne j$ we take quantum analogues of the corresponding root vectors of 
$sl_N$ multiplied by some polynomials in the diagonal generators of 
$\U_q(sl_N)$. We assume that the vector space of left-invariant 1-forms 
has 
dimension $N^2-1$ and that $X_{ij} (u^r_s)=\delta_{ir}\delta_{js}$ for
$i\ne j$. In case of the ordinary differential calculus on $SL(N)$ we
have $\sum_i \omega_{ii}=0$, so it seems to be natural to suppose that 
$\sum
q^{-2i}\omega_{ii}=0$ in the quantum case. Then all functionals of the 
corresponding quantum Lie algebra annihilate the quantum trace 
$U := \sum q^{-2i}u^i_i$.
Note that Woronowicz' 3D-calculus on $SU_q(2)$ fits into this scheme, see
Section 2 for details.\nz
This paper is organized as follows. Section 1 contains some general 
results
about left-covariant differential calculi on quantum groups which will
be needed later. In particular, we describe the construction of the 
universal
higher order differential calculus associated with a given left-covariant
first order calculus on a quantum group.\nz
In Section 2 we develop four left-covariant differential calculi 
$(\Gamma_r,\d),
r=1,2,3,4$, on the quantum group $SL_q(2)$ which satisfy the above
requirements. All four first order calculi and quantum Lie algebras give 
the
ordinary differential calculus on $SL(2)$ and the Lie algebra $sl_2$ when
$q\rightarrow 1$. However, this changes if we look at the associated 
higher
order calculi. For only one of these calculi the higher order  calculus
yields the classical calculus on $SL(2)$ in the limit $q\rightarrow1$.
As might be expected, this is Woronowicz' 3D-calculus on $SU_q(2)$ or 
more
precisely its analogue for $SL_q(2)$. For the other three calculi the
2-form $\omega_2\wedge \omega_0$ is zero and all 3-forms vanish.\nz
In Sections 3 and 4 we are concerned with left-convariant differential
calculi on the quantum group $SL_q(3)$. The functionals $X_n$ and 
$X_{ij}$ with $|i-j|=1$ are defined completely similar to the 
corresponding formulas for the 3D-calculus in Section 2. For the 
functionals 
$X_{13}$ and $X_{31}$ we use the Ansatz $X_{13}=X_{12}X_{23}-
\alpha X_{23}X_{12}$ and $X_{31}=X_{32}X_{21}-\beta X_{21}
X_{23}$ with $\alpha$ and $\beta$ complex. For arbitrary complex 
parameters
$\alpha$ and $\beta$, we obtain a first order differential calculus on 
$SL_q(3)$ which fits into the above scheme. It turns out that if 
$(\alpha,\beta)\ne (q^{-1},q)$ and $(\alpha,\beta)\ne (q,q^{-1})$, then
the 2-form $\omega_{31}\wedge \omega_{13}$ vanishes and hence 
the space
of 2-forms does not yield the corresponding space for the classical 
differential calculus on $SL(3)$ when $q\rightarrow 1$. The differential
calculi obtained in the two remaining cases $(\alpha,\beta)=(q,q^{-1})$
and $(\alpha,\beta)=(q^{-1},q)$ for the parameters $\alpha$ and 
$\beta$ are
studied in Section 4. In both cases the higher order calculi give the
ordinary higher order calculus on $SL(3)$ in the limit $q\rightarrow 1$.
The corresponding formulas show that these two calculi are
very close to the classical differential calculi on $SL(3)$ in many 
respects.\nz
In Section 5 we generalize the two differential calculi on $SL_q(3)$ from
Section 4 to $SL_q(N)$. We define two $N^2-1$ dimensional left-covariant
first order differential calculi $(\Gamma_r,\d), r=1,2$, over $SL_q(N)$
for which all quantum Lie algebra generators annihilate the quantum 
trace.
Both first order calculi give the ordinary first order calculus on 
$SL(N)$
when $q\rightarrow 1$. If $N\ge 4$, this is no longer
true for the higher order calculi.\sn
Throughout this paper $q$ is a nonzero complex number such that 
$q^2\ne 1$
and we abbreviate $\lambda:=q-q^{-1}$ and $\lambda_+:=q+q^{-1}$.\nz
We recall the definition of the quantized universal enveloping algebra
$\U_q(sl_N)$, see [4],[7]. We shall need it only for $N=2$ and $N=3$.
The algebra $\U_q(sl_{N})$ has $4(N-1)$ generators
$k_i,k^{-1}_i,e_i,f_i,i=1,\dots, N-1$, with defining relations 
\vfill\eject\noindent
$$\eqalign{
&k_ik^{-1}_i=k^{-1}k_i=1, k_ik_j=k_jk_i,k_ie_i=qe_ik_i,k_if_i=
q^{-1}f_ik_i,\cr
&e_if_j-f_je_i=\delta_{ij}\lambda^{-1}(k^2_i-k^{-2}_i),\cr
&k_ie_j=q^{-1/2}e_jk_i~{\rm and}~k_if_j=
q^{1/2} f_jk_i~{\rm if}~|i-j|=1,\cr
&e^2_ie_j-\lambda_+e_ie_je_i+e_je^2_i=f^2_if_j-\lambda_+
f_if_jf_i+f_jf^2_i=0~{\rm if}~|i-j|=1,\cr
&k_ie_j=e_jk_i, k_if_j=f_jk_i, e_ie_j=e_je_i~{\rm and}~f_if_j=
f_jf_i~{\rm if}~|i-j|\ge 2.\cr}
$$
The Hopf algebra structure of $\U_q(sl_N)$ is given by the 
comultiplication 
$\Delta$ with $\Delta (k_i)=k_i\otimes k_i,\Delta (e_i)=k_i\otimes 
e_i+e_i
\otimes k^{-1}_i,
\Delta (f_i)=k_i\otimes f_i+f_i\otimes k^{-1}_i$ and the counit 
$\varepsilon$
with $\varepsilon (k_i)=1,\varepsilon (e_i)=\varepsilon (f_i)=0$. 
There is a pairing between
the Hopf algebras $\U_q(sl_N)$ and $SL_q(N)$ such that 
$(k_i,u^n_m)=\delta_{mn}$ if
$n\ne m$ or $n=m\ne i,i+1,  (k_i,u^i_i)=
q^{1/2},(k_i,u^{i+1}_{i+1})=q^{-1/2}, (e_i,u^n_m)=\delta_{ni}
\delta_{m,i+1}$ and $(f_i,u^n_m)=
\delta_{n,i+1}\delta_{mi}$, where $u^n_m$ are the matrix entries of the
fundamental matrix of $SL_q(N)$.
\bn
{\bf 1. Left - Covariant Differential Calculi on Quantum Groups}
\mn
Our basic reference concerning differential calculi on quantum groups is 
[18].\nz
Let ${\cal A}$ be a fixed Hopf algebra with comultiplication $\Delta$, 
counit
$\varepsilon$, antipode $\kappa$ and unit element 1. Sometimes we use 
Sweedler's notation 
$\Delta^{(n)} (a) = a_{(1)} \otimes a_{(2)}\otimes\cdots\otimes 
a_{(n+1)}.$
For $a\in{\cal A}$ we put $\tilde{a} := a -\varepsilon (a)1.$\nz 
A {\it first order
differential calculus} (briefly, a FODC) over ${\cal A}$ is a pair 
$(\Gamma$,d)
of an ${\cal A}$-bimodule $\Gamma$ and a linear mapping 
d: ${\cal A}\to\Gamma$
such that d$(ab){=}$ d$a{\cdot} b {+} a{\cdot}$ d$b$ for
 $a, b\in{\cal A}$ and
$\Gamma =$ Lin $\left \{ a\cdot {\rm d} b: a, b\in{\cal A}\right \}.$ A
FODC $(\Gamma , {\rm d})$ is called {\it left covariant} if there is a 
linear mapping $\Delta_L :\Gamma \to {\cal A}\otimes \Gamma$ for which
$\Delta_L (a {\rm d}b) = \Delta (a) ({\rm id} \otimes {\rm d}) 
\Delta (b), a, b\in{\cal A}$.\nz
Suppose that $(\Gamma, {\rm d})$ is a left-covariant FODC over 
${\cal A}$. Recall that the canonical projection of $\Gamma$ into 
$\Gamma_{inv} :=\left\{ \omega\in\Gamma :\Delta_L (\omega) = 
1\otimes\omega\right\}$ is defined by $P_{{inv}} ({\rm d} a) = \kappa 
(a_{(1)}) {\rm d}a_{(2)}$, cf. [18]. We abbreviate 
$\omega(a) = P_{{inv}} ({\rm d}a).$ Then ${\cal R} := 
\left\{x\in~{\rm ker}~\varepsilon :\omega (x) = 0\right\}$ is the
{\it right ideal} of ker~$\varepsilon$  {\it associated with} 
$(\Gamma, {\rm d})$. The vector space
${\cal X} := \left\{ X\in{\cal A}' : X(1) = 0~{\rm and}~ X(x) =
 0~{\rm for}~
x\in{\cal R}\right \}$ is called the {\it quantum Lie algebra} of the
 FODC  $(\Gamma, {\rm d}).$
\bn
{\bf Lemma 1.} {\it A vector space ${\cal X}$ of linear functionals
on ${\cal A}$ is the quantum Lie algebra of a left-covariant FODC
$(\Gamma, {\rm d})$ if and only if $X(1) = 0$ and $\Delta X -\varepsilon
\otimes X \in {\cal X}\otimes {\cal A}'$ for all $X\in{\cal X}$.}
\mn
Proof. The necessity of the condition $\Delta X -\varepsilon\otimes 
X\in{\cal X}\otimes {\cal A}'$ follows at once {}{}From formula (5.20)
 in [18].
To prove the sufficiency part, let us note that the above conditions 
imply
that ${\cal R} := \left\{ x\in~{\rm ker}~\varepsilon : X(x) = 0~{\rm for}~
X\in{\cal X}\right \}$ is a right ideal of ker $\varepsilon$. {}Fr{}om the 
general
theory (cf. Theorem 1.5 in [18]) we conclude easily that ${\cal R}$ is 
the
right ideal associated with some left-covariant FODC over ${\cal A}$.\qed
\mn
The calculus $(\Gamma, {\rm d})$ is uniquely determined by ${\cal X}$ 
(because ${\cal R}$ is so) and can be described as follows. Let
$\left\{X_i: i\in I\right \}$ be a basis of the vector space ${\cal X}$ 
and
$\left\{x_i: i\in I\right\}$ a set of elements of ${\cal A}$ such that 
$X_i (x_j) = \delta_{ij}.$ Then, letting $\omega_i =\omega(x_i)$, we have
$${\rm d}a = \sum \left ( X_i * a\right ) \omega_i, a\in{\cal A}.$$
\bn
For notational simplicity we shall write $\eta\otimes \zeta$ instead of
$\eta\otimes_{{\cal A}} \zeta$, where $\eta, \zeta\in\Gamma$. We set for
$x\in{\cal R}$
$${\cal S}(x) :=\sum_{i, j}\left(X_i X_j\right) (x)~\omega_i
\otimes\omega_j.$$
\bn
Some properties of the mapping 
${\cal S} : {\cal R}\to\Gamma\otimes_{{\cal A}}\Gamma$ are collected in 
the following
\bn
{\bf Lemma 2.} {\it Let $P_{{inv}}$ denote the canonial 
projection of the left-covariant bimodule $\Gamma
\otimes_{{\cal A}}\Gamma$
into  $(\Gamma\otimes_{{\cal A}}\Gamma)_{{inv}}$. For $x\in{\cal R}$ and
$a\in{\cal A}$, we have:}
\mn
(i) $\S1(x) a = a_{(1)} \S1 (x a_{(2)}).$\nz
(ii) $\S1(x a) = \kappa (a_{(1)}) \S1(x) a_{(2)} = P_{{inv}} 
(\S1 (x) a).$
\par
\bn
{\bf Proof.}~We show the first equality of (ii). Recall that 
$\Delta X_i=\varepsilon \otimes X_i+X_k\otimes f_i^k, i\in I$, by 
formula (5.20) in [18], where
$f_i^k$ are functionals on ${\cal{A}}$ such that $\omega_k a=
(f_i^k\ast a)\omega_i$, $a\in {\cal{A}}$. Hence we get 
$$\eqalign{
\S1 (xa)&=\sum_{i,j} (X_iX_j)(xa)~\omega_i\otimes\omega_j=
\sum_{i,j} (\Delta X_i)(\Delta
X_j)(x\otimes a)~\omega_i\otimes\omega_j\cr
&=\sum_{i,j,k,l} (X_k X_l)(x) f^k_i (a_{(1)})
f^l_j (a_{(2)})~\omega_i\otimes \omega_j~.\cr}
$$
On the other hand, we have 
$$\eqalign{
\kappa (a_{(1)})\S1 (x) a_{(2)}&=\sum_{k,l}\kappa(a_{(1)})
(X_k X_l)(x)\omega_k\otimes\omega_l a_{(2)}\cr
&=\sum_{k,l}\kappa (a_{(1)}) (X_kX_l)(x)
f_i^k\ast (f_j^l\ast a_{(2)})\omega_i\otimes\omega_j\cr
&=\sum_{i,j,k,l}\kappa (a_{(1)}) a_{(2)}
f^k_i (a_{(3)}) f^l_j (a_{(4)}) (X_k X_l)(x)\omega_i\otimes\omega_j~.\cr}
$$
By the Hopf algebra axioms both expressions are equal.\nz
The preceding calculation yields $a_{(1)}\S1 (xa_{(2)})=a_{(1)}
\kappa (a_{(2)})\S1 (x)a_{(3)}=\S1 (x) a$ which proves (i). 
Applying $P_{inv}$ to (i) and
using the fact that $P_{inv} (a\eta)=\varepsilon (a)P_{inv}(\eta)$ 
we get the second equality of (ii).~\qed
\bn
A {\it differential calculus} over ${\cal A}$ is a pair $(\Gamma^\land, 
{\rm d})$ of a graded algebra $\Gamma^\land = \oplus_{n=0}^\infty 
\Gamma_n^\land$
with product $\land : \Gamma_n^\land\times\Gamma_m^\land\to
\Gamma_{n+m}^\land$
and a linear mapping d$ : \Gamma^\land\to\Gamma^\land$ of degree one 
such that
$\d^2 = 0,  \d:\Gamma_n^\land \to\Gamma_{n+1}^\land, \d(\eta\land\zeta) =
\d\eta\land\zeta + (-1)^n \eta\land$ d$\zeta$ for $\eta\in\Gamma_n^\land, 
\zeta\in\Gamma_m^\land, \Gamma_0^\land = {\cal A}$ and 
$\Gamma_n^\land = {\rm Lin} \left \{ a~ {\rm d}a_1 \wedge\ldots \wedge
{\rm d}a_n: 
a, a_1,\ldots, a_n\in{\cal A}\right \}$ for $n\in\bbbn.$ The definition
 of
left-covariance of a differential calculus is similar to the case of 
first order calculi. Sometimes we simply write $\eta\zeta$ for $\eta
\wedge \zeta$.\nz
For each left-covariant FODC $(\Gamma, {\rm d}_1)$ over ${\cal A}$ there
exists a unique (up to isomorphism) universal left-covariant differential
calculus $(\Gamma^\land, {\rm d})$ over ${\cal A}$ such that 
$\Gamma_1^\land = \Gamma$ and d$\lceil {\cal A} =$ d$_1 .$ We briefly
describe
the construction of $(\Gamma^\land, {\rm d}$).\nz\hfill\eject\noindent
Let $\Omega = \oplus_{n=0}^\infty \Omega^n$ be the universal differential
envelope of the algebra ${\cal A}$, see e.g. [2] or [3]. We have 
$\Omega^0 = {\cal A}$ and $\Omega^n = {\cal A}\otimes ({\rm ker}~
\varepsilon)^{\otimes n}$ by identifying $a\otimes a_1\otimes 
\ldots\otimes a_n$
and $a\d a_1\ldots\d a_n.$ The differential d of $\Omega$ is given by 
$\d(a\d a_1\ldots \d a_n) =\d a \d a_1\ldots\d a_n$. Clearly, 
$(\Omega, {\rm d})$ is a differential calculus over the Hopf algebra 
${\cal A}$.
Let ${\cal J} ({\cal R}) :=\Omega~\omega ({\cal R})\Omega +
 \Omega~$d$\omega
({\cal R})\Omega$ be the differential ideal generated by the set
$\omega({\cal R})$. Here ${\cal R}$ is the right ideal of
 ker $\varepsilon$
associated with the given FODC $(\Gamma, {\rm d}_1)$. We have
${\cal J} ({\cal R}) =\sum_n {\cal J}_n ({\cal R})$, where
${\cal J}_n ({\cal R}) := {\cal J} ({\cal R})\cap \Omega^n.$ Obviously, 
the
quotient algebra $\Omega / {\cal J} ({\cal R}) = \sum_n\Omega^n /
{\cal J}_n ({\cal R})$ endowed with the quotient map of d is also 
differential
calculus over ${\cal A}$. By formula (1.23) in [18], ${\cal J}_1 
({\cal R}) =
{\cal A}\omega ({\cal R})$ coincides with the submodule ${\cal N}$ 
occuring
in Theorem 1.5 of [18]. Therefore, the first order calculus 
$(\Omega_1/{\cal J}_1
({\cal R}),\d)$ of $(\Omega / {\cal J} ({\cal R}),\d)$ is isomorphic to
$(\Gamma, {\rm d}_1)$. Moreover, it is not difficult to verify that 
$(\Omega / {\cal J} ({\cal R}),$ d) is left-covariant. {}From the preceding
construction it is clear that $(\Omega / {\cal J} ({\cal R}),$ d) has the 
following universal property: If $(\Omega', {\rm d}')$ is another 
differential
calculus over ${\cal A}$ such that  $(\Omega_1^{'}, {\rm d}')$ is 
isomorphic
to the FODC $(\Gamma, {\rm d}_1)$, then  $(\Omega', {\rm d}')$ is
(isomorphic to) a quotient of $(\Omega / {\cal J} ({\cal R}),$ d) by some
differential ideal.\nz 
In order to obtain a more explicit description of the calculus
$(\Omega / {\cal J},\d)$, we shall use a construction of the 
differential 
envelope $\Omega = \oplus_n \Omega^n$ of the Hopf algebra
${\cal A}$ developed in [14]. More details and proofs of all unproven 
assertions
in the following discussion can be found in [14].\nz
Let $\Omega^0 := {\cal A}.$ For
$n\in\bbbn$, we set $\Omega^n := {\cal A} \otimes (\ker 
\varepsilon)^{\otimes n}$ and we write $a \omega (a_1)\ldots \omega(a_n)$
instead of $a\otimes a_1\otimes\ldots\otimes a_n.$ We put 
$\omega (\lambda 1 + a) =\omega(a)$ for $\lambda\in\bbbc, a\in\ker 
\varepsilon.$
The product of $\Omega =\oplus_{n=0}^\infty\Omega^n$ and the 
differential d are defined by
\bn
$$\eqalignno{
&a\omega(a_1)\ldots\omega (a_n) b~\omega (b_1)\ldots\omega (b_m) :=
\quad\cr
&ab_{(1)}~\omega(a_1 b_{(2)}){\ldots}~\omega (a_{n{-}2}
 b_{(n{-}1)})~\omega 
(a_{n{-}1} b_{(n)})~\omega (a_n b_{(n{+}1)})~\omega (b_1){\ldots}~
\omega (b_m)\cr}
$$
$$\eqalign{
{\rm and}~\d (a\omega (a_1)\ldots\omega (a_n)) {:=}~ &a_{(1)}\omega 
(a_{(2)})~\omega (a_1)
\dots\omega (a_n){+}\cr
&{\sum\limits^n_{i{=}1}} (-1)^ia\omega(a_1){\dots}\omega(a_{i{-}1})\omega 
(a_{i,(1)})\omega (a_{i,(2)})\omega(a_{i{+}1}){\dots}\omega(a_n),\cr}
$$
where $a_1,{\dots}, a_n,b_1, {\dots}, b_m \in \ker\varepsilon$ and 
$a,b\in {\cal{A}}$.~
(In order to motivate these formulas, we recall that for any 
left-covariant differential
calculus over ${\cal A}$ we have $\d a=a_{(1)}\omega (a_{(2)}),\d
\omega (a)=-
\omega(a_{(1)})\wedge \omega (a_{(2)})$ and $\omega (b)c=c_{(1)}\omega 
(bc_{(2)})$
for $a,c\in{\cal A}$ and $b\in\ker\varepsilon$.) It can be shown that the
 pair
$(\Omega,\d)$ endowed with the above definitions becomes a 
left-covariant 
differential calculus over ${\cal A}$ which is (isomorphic to) the
differential envelope of ${\cal A}$. For $a\in\R$, the element $\omega (a)$ of 
$\Omega$ is equal to $\kappa (a_{(1)})\d a_{(2)}$ which justifies to 
use the notation $\omega (a)$. Obviously, the kernel of the map $\omega:
{\cal A}\rightarrow\Omega$ is $\bbbc\cdot 1$. \nz
Let $\R,\{X_i\},\{x_i\},\{\omega_i\}$ and $\S1$ be as
defined above for the left-covariant FODC $(\Gamma,\d_1)$. We put
$$
\S1_u(x):=\sum\limits_{i,j} (X_iX_j)(x)\omega_i\omega_j~{\rm for}~
x\in\R~,
$$
where the product $\omega_i\omega_j$ is taken in the algebra $\Omega$. For
$a\in{\cal {A}}$, the element $\tilde{a}-\sum_i X_i (a)\tilde{x}_i$ is 
annihilated by $\X$ and by $\varepsilon$, so it belongs to $\R$ and hence
$\omega (a)-\sum_i X_i (a)\omega_i\in\omega (\R)$. Therefore, since
$\d \omega (x)=-\omega (x_{(1)})\omega (x_{(2)})$, we obtain
$$\eqalign{
&\d\omega (x)+\S1_u(x)=\cr
&-\left( \omega (x_{(1)})-\sum\nolimits_iX_i(x_{(1)})\omega_i\right)~
\omega(x_{(2)})+\sum\nolimits_iX_i
(x_{(1)})\omega_i\left( \sum\nolimits_jX_j(x_{(2)})\omega_j-
\omega (x_{(2)})\right)\cr
&\in\omega (\R)\Omega^1+\Omega^1\omega (\R)~{\rm for}~x\in\R~.\cr}
$$
Hence the differential ideal ${\cal J}(\R)$ is generated by the sets $\omega
(\R)$ and $\S1_u(\R).$\nz
Next we define the exterior algebra for the FODC $(\Gamma,\d_1)$. Let
$\Gamma^{\otimes 0}:={\cal A},\Gamma^{\otimes n}:=\break
\Gamma\otimes_{\cal A}{\dots}\otimes_{\cal A}
\Gamma$ ($n$ times) for $n\in\bbbn$ and let $\Gamma^\otimes:=
{\oplus^\infty_{n=0}}\Gamma^{\otimes n}$ be the tensor algebra
of $\Gamma$ over ${\cal A}$. We denote by $S={\oplus^\infty_{n=2}}S_n$ the
two-sided ideal of the algebra $\Gamma^\otimes$ generated by the set
$\S1(\R)$. The quotient algebra $\Gamma^\wedge=\Gamma^\otimes/S$ is called
the {\it exterior algebra} over ${\cal A}$ {\it for the} FODC $(\Gamma,\d_1)$.
Clearly, $\Gamma^\wedge$ is also a graded algebra $\Gamma^\wedge={\oplus^\infty_{n=0}}
\Gamma^\wedge_n$ with $\Gamma^\wedge_0={\cal A},\Gamma^\wedge_1=\Gamma$ and
$\Gamma^\wedge_n=\Gamma^{\otimes n}/S_n$ for $n\in \bbbn$. The product of
$\Gamma^\wedge$ is denoted by $\wedge$.\nz
For $x\in\ker\varepsilon,$ let $[x]$ denote the coset $x+\R$. Let $\pi$
be the product of the mapping $\pi_1:\Omega\rightarrow\Gamma^\otimes$
defined by $\pi_1(a)=a, \pi_1(a\omega (a_1)\dots\omega(a_n))=a\omega([a_1])\dots
\omega ([a_n])$ and the quotient map {}From $\Gamma^\otimes$ onto $\Gamma^\wedge$.
Then $\pi$  is an algebra homomorphism of $\Omega$ onto $\Gamma^\wedge$.
It can be shown that the kernel of $\pi$ is the two-sided ideal in
$\Omega$ generated by $\omega(\R)$ and $\S1_u(\R)$. Therefore, by the
paragraph before last, $\ker\pi={\cal J} (\R)$. Hence the quotient algebra
$\Omega/{\cal J} (\R)$ and the exterior algebra $\Gamma^\wedge$ are isomorphic.
We define $\d\pi (a):=\d a$ for $a\in\Gamma$. Then $(\Gamma^1,\d)$ is a
left-covariant differential calculus over ${\cal A}$ which is isomorphic to 
the calculus $(\Omega/{\cal J} (\R),\d)$. By construction, the first order
calculus $(\Gamma^\wedge_1,\d)$ of $(\Gamma^\wedge, \d)$ is the given FODC
$(\Gamma,\d_1)$. We call $(\Gamma^\wedge, \d)$ (or likewise $(\Omega/{\cal J}
(\R),\d)$ the {\it universal differential calculus associated with the} FODC
$(\Gamma,\d_1)$.\nz
By definition, $S_2={\cal A}\S1 (\R) {\cal A}$. {}From Lemma 2, (i), we see
that $S_2={\cal A}\S1 (\R)$ is an ${\cal A}$-subbimodule of $\Gamma\otimes_{\cal A}\Gamma$ and
hence a left-covariant bimodule over ${\cal A}$. Therefore, each basis $(\zeta_n)$
of the vector space $\S1 (\R)$ of symmetric elements is a free left module
basis for $S_2$, i.e. any element $\zeta$ of $S_2$ can be written as
$\zeta=\sum_n a_n\zeta_n$ with elements $a_n\in{\cal A}$ uniquely determined by
$\zeta.$
\bn
\noindent
{\bf 2. Left-covariant differential calculi on ${\bf SL_q(2)}$}
\mn
In this section we denote the matrix entries $u^1_1,u^1_2,u^2_1,u^2_2$ for the
quantum group $SL_q(2)$ by $a,b,c,d$, respectively, and the generators of
$\U_q(sl_2)$ by $k,k^{-1}$, $e,f$.
Our aim is to study left-covariant differential calculi $(\Gamma,\d)$ on
${\cal A}:=SL_q(2)$ of the form
$$\eqalignno{
&\d x=\sum\limits^2_{i=0}(X_i\ast x) ~\omega_i,~x\in SL_q(2)~, &(2.1)\cr}
$$
where $\omega_0, \omega_1$ and $\omega_2$ are left-invariant 1-forms and
$X_0,X_1$ and  $X_2$ are linear functionals {}From $\U_q(sl(2))$
satisfying
$$\eqalignno{
&X_1(a){=}X_0(b){=}X_2(c){=}1,X_1(b){=}X_1(c){=}X_0(a){=}X_0(c){=}
X_2(a){=}X_2(b){=}0\hskip1cm &(2.2)\cr
{\rm and}\hskip0.2cm &X_i(1){=} X_i(q^{-2} a+q^{-4}d){=}
0~{\rm for}~i{=}0,1,2~.\hskip01cm&(2.3)\cr}
$$
{}From the pairing between $\U_q(sl_2)$ and $SL_q(2)$ it follows that arbitrary linear
functionals $X_i\in\u, i=0,1,2$, satisfying (2.2) and (2.3) can be
written as $X_1=efp_{11} (k)+p_{12}(k)+X'_1, X_0=e p_0(k)+X'_0$
and $X_2=fp_2(k)+X'_2$, where $p_{11}, p_{12}, p_0, p_2$ are Laurent
polynomials in $k$ and $X'_1,X'_0,X'_2$ annihilate all four matrix
entries $a,b,c,d$. This suggests to consider the following Ansatz
$$\eqalignno{
&X_1=ef p_{11} (k)+p_{12} (k), X_0=e p_0(k), X_2 =
fp_2(k)&(2.4)\cr}
$$
with polynomials $p_{11},p_{12},p_0$ and $p_2$ in $k$ and $k^{-1}$.
\bn
{\bf Theorem 1.}~{\it There are precisely} 4 {\it non-isomorphic three dimensional
left-covariant differential claculi $(\Gamma_r,\d), r=1,2,3,4$,
satisfying (2.1)-(2.3) obtained by the Ansatz (2.4). The right ideal
$\R_r$ of $\ker\varepsilon$ associated with the calculus $(\Gamma_r, \d)$ is
generated by the six elements
$$\eqalign{
&a+q^{-2}d-(1+q^{-2})1,~b^2,~c^2,~bc,~(a-\gamma_{0r})
b,~(a-\gamma_{2r})c,~\cr}
$$
where $\gamma_{0r}$ and $\gamma_{2r}$ are the coefficients given by
$\gamma_{01}=\gamma_{02}=1,\gamma_{03}=\gamma_{04}=q^{-2},\gamma_{21}=
\gamma_{23}=1$ and $\gamma_{22}=\gamma_{24}=q^{-2}$.}\sn
\bn
{\bf Proof.} We first suppose that $(\Gamma,\d)$ is a differential
calculus such that (2.1)-(2.4) are fulfilled. Let $\R$ be its associated
right ideal of $\ker\varepsilon$. {}From the formulas (2.1)-(2.4) and {}From the
pairing between $\U_q(sl_2)$ and $SL_q(2)$ we compute easily that $b^2$
and $c^2$ are in $\R$ and that there are complex numbers $\gamma_1,
\gamma_0,\gamma_2$ such that $bc-\gamma_1(a-1), ab-\gamma_0b$ and
$ac-\gamma_2c$ are in $\R$. As usual, we shall write $x\equiv y$ if
$x-y\in\R$. By (2.3), $a-1\equiv -q^{-2}(d-1)$. Thus
$$\eqalignno{
0&\equiv (bc-\gamma_1 (a-1))(a-1)=q^{-2}abc-bc-\gamma_1
(a-1)^2\cr
&=\gamma_0q^{-2} bc-bc-\gamma_1(a-1)^2=(\gamma_0 q^{-2}-1)
\gamma_1(a-1)-\gamma_1(a-1)^2\cr
&=-\gamma_1(a-1)(a-\gamma_0q^{-2}1)=\gamma_1q^{-2} (d-1)
(a-\gamma_0q^{-2}1)\cr
&=\gamma_1 q^{-2}(da-a-\gamma_0q^{-2}(d-1))=\gamma_1q^{-2} (q^{-1}
bc+1-a-\gamma_0 q^{-2} (d-1))\cr
&\equiv \gamma_1 q^{-2} (\gamma_1q^{-1} (a-1)+1-a+\gamma_0
(a-1))=\gamma_1q^{-2} (\gamma_1q^{-1}+\gamma_0-1)(a-1)~\cr
{\rm and}\hskip0.6cm&&\cr
0&\equiv (bc-\gamma_1 (a-1))b=b^2c-\gamma_1 ab+\gamma_1b\equiv \gamma_1
(1-\gamma_0)b~.\cr}
$$
Since $a-1, b$ and $c$ are  not in $\R$ by (2.2), we get $\gamma_1
(\gamma_1q^{-1}+\gamma_0-1)=0$ and $\gamma_1(1-\gamma_0)=0$, so that
$\gamma_1=0.$ This implies that $da-1\equiv 0$ and hence $q^2(a-1)
(a-q^{-2})=(d-1)(q^{-2}-a)\equiv -1 +a+q^{-2} (d-1)\equiv 0$,
so $0\equiv (a-1)(a-q^{-2})b=qaba-(1+q^{-2})ab+q^{-2}b=(\gamma_0-1)
(\gamma_0-q^{-2})b$ which yields $(\gamma_0-1)(\gamma_0-q^{-2})=0$. Similarly
we obtain $(\gamma_2-1)(\gamma_2-q^{-2})=0$. The two latter equations imply
that $\R$ contains one of the right ideals $\R_r,r=1,2,3,4$. {}From (2.3),
$\omega_0,\omega_1$ and $\omega_2$ are linearly independent. Hence
we have ${\rm codim}~\R=\dim\Gamma_{inv}\ge 3$.
Since ${\rm codim}~ \R_r\le 3$ by the definition of $\R_r$, we
conclude that $\R =\R_r$.\sn
To complete the proof of Theorem 1, we have to construct the first order
differential calculi $(\Gamma_r,\d)$ having the desired properties. For
this let $\X_r$ denote the linear span of functionals $X_1,X_0,X_2
\in\U_q(sl_2)$, where
$$\eqalign{
&X_0 := q^{-1/2} ek^{-1}~{\rm for}~ r=1,2,~X_0 := q^{-5/2}
 ek^{-5}~{\rm for}~ r=3,4~,\hskip2.6cm\cr
&X_2 := q^{1/2} fk^{-1}~{\rm for}~ r=1,3, ~X_2 := q^{5/2} fk^{-5}~
{\rm for}~ r=2,4\cr
{\rm and}\hskip0.5cm~ &X_1 := q\lambda^{-1} 
(\varepsilon-k^{-4})~{\rm for}~r=1,2,3,4~.
\cr}
$$
{}From these definitions and {}From the comultiplication in $\U_q(sl_2)$ 
we obtain
$$\eqalign{
&\Delta X_j=\varepsilon\otimes X_j+X_j\otimes k^{-2}~ {\rm for}~
r=1,2, j=0~~{\rm and}~ r=1, j=2,\cr
&\Delta X_j =\varepsilon \otimes X_j+X_j\otimes
k^{-6} +(q^{-2}-1)X_1\otimes X_j ~{\rm for}~r=3,4, j=0~
{\rm  and}~r=2,4, j=2,\cr
&\Delta X_1=\varepsilon\otimes X_1+X_1\otimes k^{-4}~{\rm for}~
r=1,2,3,4~.\cr}
$$
This shows that $\Delta X-\varepsilon\otimes X\in\X_r\otimes {\cal {A}}'$ for
all $X\in\X_r$. Therefore, by Lemma 1, each vector space $\X_r$ defines a left-
covariant first order differential calculus $(\Gamma_r,\d)$ over $SL_q(2)$.
Let $\R'_r$ denote the right ideal of $\ker\varepsilon$ associated with
the calculus $(\Gamma_r,\d)$. One verifies that the functionals $X_0,
X_1,X_2$ for the calculus $(\Gamma_r,\d)$ have the properties (2.2)-(2.3)
and that they annihilate the six generators of the right ideal $\R_r$. Hence
$\R_r'\subseteq\R_r$. Since ${\rm codim}~ \R'_r={\rm dim}~ (\Gamma_r)_{inv}=3$
and ${\rm codim}~\R_r\le 3$, we have $\R'_r=\R_r$. This completes the proof
of Theorem 1.~\qed\sn
We now describe the structure of the four differential calculi $(\Gamma_r,\d)$,
$r=1,2,3,4$, more in detail. By the general theory [18], the above formulas
for the comultiplication of $X_j$ lead to the following commutation rules
between matrix entries and 1-forms:
$$
\eqalign{
&\omega_j a{=}q^{-1} a\omega_j, \omega_j b{=} q b \omega_j, \omega_j c{=}
q^{-1} c \omega_j, \omega_j d{=}q d \omega_j~{\rm for}~ r{=}1,2, j{=}0~
{\rm and}~ r{=}1,3,j{=}2,\cr
&\omega_j {a{=}q^{-3}} a\omega_j, \omega_j {b{=}q^3} b \omega_j, \omega_j
{c{=}q^{-3}} c \omega_j, \omega_j {d{=}q^3} d \omega_j~{\rm for}~ {r{=}3,4},
{j{=}0}~{\rm and}~ {r{=}2,4,j{=}2},\cr
&\omega_1 a= q^{-2} a\omega_1, \omega_1 c = q^{-2} c\omega_1~
{\rm for}~ r=1, 3,\cr
&\omega_1 b= q^{2} b\omega_1, \omega_1 d= q^2 d \omega_{1}~{\rm for}~ r=1,2,\cr
&\omega_1 a= q^{-2} a\omega_1 + (q^{-2}-1)b\omega_2, \omega_1 c= q^{-2}
c \omega_1 +(q^{-2}-1)d\omega_2~{\rm for}~r=2,4,\cr
&\omega_1 b= q^{2} b\omega_1 + (q^{-2}-1) a\omega_0, \omega_1 d= q^{2}
d \omega_1+ (q^{-2}-1)c \omega_0~{\rm for}~ r=3,4~.\cr}
$$
Recall that $\omega(a)=P_{inv}(\d a)=\kappa (a_{(1)}) \d a_{(2)}$
if $\Delta (a)=a_{(1)}\otimes a_{(2)}$. {}From the formulas (2.2) we compute
that for all four calculi
$$\eqalign{
&\omega_1= \omega (a), \omega_0= \omega (b), \omega_2= \omega (c),\cr
{\rm and}~\hskip1cm &\d a=b \omega_0+a\omega_1, \d b=a\omega_2-q^2
b\omega_1, \d c= c \omega_1 + \d\omega_0, \d d=-q^2d\omega_1+c\omega_2.\qquad\cr}
$$
According to the general theory [18], the calculus $(\Gamma_r,\d)$ is a $\ast$-
calculus for an algebra involution $x\rightarrow x^\ast$ on $SL_q(2)$ if and only
if $\kappa (x)^\ast\in\R_r$ for all $x\in\R_r$. Obviously, it suffices to check
this condition for the six generators of the right ideal $\R_r$.
The four calculi $(\Gamma_r,\d), r=1,2,3,4$, are $\ast$-calculi for the Hopf
$\ast$-algebra $SL_q(2,\bbbr), |q|=1$, while only $(\Gamma_1,\d)$ and
$(\Gamma_4, \d)$ are $\ast$-calculi for the real forms $SU_q(2)$ and $SU_q (1,1),
q\in\bbbr$, of the quantum group $SL_q(2)$. However, for the Hopf $\ast$-
algebras $SU_q(2)$ and $SU_q(1,1)$ we have $\kappa (x)^\ast\in\R_2$
for $x\in \R_3$ and $\kappa (x)^\ast\in \R_3$ for $x\in\R_2$.
Moreover, we have $\varphi (\R_2)=\R_3$ and $\varphi (\R_3)=\R_2$, where
$\varphi$ denotes the algebra automorphism of $SL_q(2)$ which fixes $a$ and
$d$ and interchanges $b$ and $c$.\nz
Next we consider the commutation rules between the generators $X_0,X_1,
X_2$ of the quantum Lie algebra $\X_r$ of the calculus $(\Gamma_r,\d)$.
We have
$$\eqalign{
&q^2 X_1X_0 - q^{-2} X_0X_1 =(1+q^2)X_0~~~{\rm for}~ r=1,2,3,4,\cr
&q^2 X_2X_1 - q^{-2} X_1X_2 =(1+q^2)X_2~~~{\rm for}~ r=1,2,3,4,\cr
&q X_2X_0 - q^{-1} X_0X_2 =-q^{-1}X_1~~~{\rm for}~ r=1,\cr
&q^3X_2X_0 - q^{-3} X_0X_2 =-q^{-1}X_1+q^{-2} \lambda X^2_1~~~
{\rm for}~ r=2,3,\cr
&q^5 X_2X_0 - q^{-5} X_0X_2 =-q^{-1}X_1+2q^{-2}\lambda X^2_1-q^{-3}
\lambda^2X^3_1~~~{\rm for}~ r=4,\cr}
$$
What about the classical limits $q\rightarrow 1$ of the calculi $(\Gamma_r,\d)$?
If we keep the basis $\{\omega_0,\omega_1,\omega_2\}$ of $(\Gamma_r)_{inv}$
fixed, all above equations make sense in the limit $q\rightarrow 1$ and we
obtain the classical first order differential calculus on the Lie group
$SL(2)$. That is, all four calculi $(\Gamma_r,\d), r=1,2,3,4$, can be
considered as deformations of the classical first order differential calculus
on $SL(2)$. In particular, the preceding equations for $ r=1, r=2,3$ and $
r=4$ define three deformations of the Lie algebra $sl_2$. Note that the
quadratic and cubic terms of $X_1$ in the two last equations vanish in the
limit $q\rightarrow 1$. It might be worth to mention that the quantum Lie
algebras $\X_2$ and $\X_3$ are isomorphic (because the commutation rules
of the generators $X_0,X_1,X_2$ for $r=2$ and $r=3$ are the same),
but the right ideals $\R_2$ and $\R_3$ are different and hence the
differential calculi $(\Gamma_2,\d)$ and $(\Gamma_3, \d)$ are not ismorphic.\nz
Clearly, for the quantum group $SU_q(2)$ the differential calculus $(\Gamma_1,
\d)$ is nothing but the 3D-calculus discovered by S.L. Woronowicz [17], because
the right ideal $\R_1$ coincides with the right ideal of the 3D-calculus,
cf. formula (2.27) in [17]. (The slight differences between some of our
formulas stated above and the corresponding formulas in [17] stem {}From the
fact that we assumed $X_2(c)=1$ by (2.2), while $X_2(c)=-q$
by formula (2.3) in [17].) \nz
Now we turn to the higher order differential calculi.\bn
{\bf Lemma 3.} {\it If $q^{12}\ne 1$, then $\omega_0\otimes\omega_2\in
\S1(\R_r)$ for $r=2,3,4$.}\bn
{\bf Proof.} Obviously, $a^2-(1+q^{-2})a+q^{-2}1\in\R_r$. Easy
computations  yield
$$\eqalign{
&{\S1 (a^2{-}(1{+}q^{-2})a{+}q^{-2}1){=}(1{+}q^{-2})(q^{-2}\omega_1
{\otimes} \omega_1 {+} (\gamma_{0r}\gamma_{2r}{-}1)\omega_0{\otimes}\omega_2)}\quad\cr
{\rm and}\hskip1.9cm&\cr
&{\S1 (a+q^{-2} d-(1+q^{-2})1)
{=}(1{+}q^2)\omega_1{\otimes}\omega_1{+}\omega_0{\otimes}
\omega_2{+}q^{-2}\omega_2 {\otimes} \omega_0~}\hskip0.5cm\cr}$$
for $r=2,3,4,$ so that
$$\eqalignno{
&q^6\omega_0{\otimes}\omega_2{+}\omega_2{\otimes}\omega_0\in \S1 (\R_4)~{\rm and}~
(q^6+q^4{-}1)\omega_0{\otimes}\omega_2{+}\omega_2{\otimes}\omega_0
\in \S1(\R_2){\cap} \S1(\R_3).~ &(2.5)\cr}
$$
Let $r=2$. Then we have $ab-b\in \R_2$ and $\S1 (ab-b)=\omega_0\otimes
\omega_1+q^{-2}\omega_1\otimes\omega_0~.$
Using Lemma 2, (ii), we compute
$$\eqalignno{
&\S1 ((ab-b)c)=P_{inv}(\S1 (ab-b)c)=P_{inv}(\omega_0\otimes
\omega_1 c+q^{-2}\omega_1\otimes\omega_0 c)&\cr
&=P_{inv}(q^{-3}c\omega_0\otimes\omega_1+q(q^{-2}-1)d\omega_0
\otimes\omega_2+q^{-5} c\omega_1\otimes\omega_0+q^{-3}
(q^{-2}-1)d\omega_2\otimes\omega_0)&\cr
&=(q^{-2}-1)(q\omega_0\otimes\omega_2+q^{-3}\omega_2\otimes
\omega_0)\in {\cal{S}}(\R_2)~.&(2.6)\cr}
$$
Since $q^6\ne 1$, (2.5) and (2.6) imply that $\omega_0\otimes\omega_2\in
\S1(\R_2)$. Interchanging the role of $b$ and $c$ we get
$\omega_0\otimes\omega_2\in \S1(\R_3)$. For $r=4$, we have
$(ac-q^{-2} c)b\in\R_4$ and
$$\eqalignno{
&\S1((ac-q^{-2} c)b)=P_{inv} ((\omega_1\otimes\omega_2+q^{-4}
\omega_2\otimes\omega_1)b)&\cr
&P_{inv} (q^5 b\omega_1\otimes\omega_2+q^3(q^{-2}-1)a\omega_0
\otimes \omega_2+qb\omega_2\otimes\omega_1+q^{-7}
(q^{-2}-1)a\omega_2\otimes \omega_0)\cr
&=(q^{-2}-1)(q^3\omega_0\otimes\omega_2+q^{-7}\omega_2
\otimes \omega_0)~.&(2.7)\cr}
$$
{}From (2.5) and (2.7) we obtain $\omega_0\otimes\omega_2\in
\S1(\R_4)$, because we assumed that $q^{12}\ne 1$.~\qed\sn
Therefore, in contrast to the classical case, the 2-form $\omega_0\wedge
\omega_2$ is zero for all three calculi $(\Gamma_r,\d),r=2,3,4$. In particular,
this implies that all 3-forms vanish and that the differential of
$\omega_1$ is zero. Hence the higher order calculi of $(\Gamma_r,\d), r=
2,3,4$, do {\it not} give the ordinary differential calculus on $SL(2)$ when
$q\rightarrow 1$. Recall that the calculus $(\Gamma_1,\d)$ is 3D-calculus of
S.L. Woronowicz [17]. As shown in [17], the higher order calculus of the
3D-calculus yields the ordinary differential calculus on $SU(2)$ (and on
$SL(2)$) in the limit $q\rightarrow 1$. Thus the considerations
in this section emphasize the distinguished role of Woronowicz' 3D-calculus
on $SU_q(2)$ resp. of the calculus $(\Gamma_1,\d)$ on $SL_q(2)$ among all three
dimensional left-covariant differential calculi on $SU_q(2)$ resp. $SL_q(2)$.
\bn
{\bf Remark.} The functionals $\tilde{X}_0 :=q^{1/2} ek, \tilde{X}_2:=q^{-1/2}
fk$ and $\tilde{X}_1 :=q \lambda^{-1}(\varepsilon -k^4)$ also satisfy the commutation
relations of the quantum Lie algebra ${\cal X}_1$. This presentation has been
found by A. Sudbery [16], see e.g. [10]. The functionals $\{\tilde{X}_0,
\tilde{X}_1,\tilde{X}_2\}$
define another left covariant FODC $(\tilde{\Gamma},\d)$ on $SL_q(2)$. Since
$\tilde{X}_1$ does not annihilate the quantum trace $q^{-2}a+q^{-4}d$, the
right ideal of $\ker \varepsilon$ associated with $(\tilde{\Gamma},\d)$ is
different {}From $\R_1$, so that the first order calculi $(\tilde{\Gamma},\d)$
and $(\Gamma_1,\d)$ are not isomorphic.
\bn
{\bf 3. Left-covariant differential calculi on ${\bf SL_q(3)}$}
\mn
In this and the following section we consider left-covariant differential
calculi $(\Gamma,\d)$ over ${\cal{A}}=SL_q(3)$ of the form
$$\eqalignno{
&\d x =\sum\limits^{3}_{i,j=1;i\ne j} (X_{ij}\ast x)~\omega_{ij}+\sum\limits^2_
{n=1} (X_n\ast x)~\omega_n,x\in SL_q(3)~. &(3.1)\cr}
$$
Here $\omega_{ij}$ and $\omega_n$ are left-invariant 1-forms and $X_{ij}$
and $X_n$ are linear functionals of $\U_q(sl_3)$ such that
$$\eqalignno{
&X_{ij}(1)=0~{\rm and}~X_{ij} (u^r_s)= \delta_{ir}\delta_{js} ~{\rm for}~
i\ne j,\cr
&X_n (u^r_s)=0~{\rm for}~ r\ne s~{\rm and}~ X_n(1)=X_n (U)=0~.
&(3.2)\cr}
$$
Recall that $U:={\sum\limits^{3}_{i=1}}q^{-2i}u^i_i$ is the quantum trace.
We define
$$\eqalignno{
&X_i:= q\lambda^{-1} (\varepsilon-k_i^{-4}),~X_{i,i+1}:=q^{-1/2} e_ik_i^{-1}~
{\rm and}~X_{i+1,i} := q^{1/2} f_ik_i^{-1}~{\rm for}~i=1,2~.\hskip0.4cm&(3.3)\cr}
$$
Let $\alpha$ and $\beta$ be complex numbers. We set
$$
X_{13}=X_{12}X_{23}-\alpha X_{23}X_{12}~{\rm and}~X_{31}=
X_{32}X_{21}-\beta X_{21}X_{32}~.
$$
Then all linear functionals $X_{ij}$ and $X_n$ satisfy the conditions
(3.2). Let $\X$ denote the vector space generated by these functionals. We
compute
$$\eqalign{
&\Delta X_n{=}\varepsilon\otimes X_n{+}X_n\otimes k^{-4}_n,~~ \Delta
X_{ij} {=}\varepsilon\otimes X_{ij}{+}X_{ij}\otimes k^{-2}_j~{\rm for}~
|i-j|{=}1,\cr
&\Delta X_{13}{=}\varepsilon\otimes X_{13}{+}X_{13}\otimes (k_1k_2)^{-2}
{+}(q{-}\alpha) X_{12}\otimes X_{23}k^{-2}_1{+}(1{-}\alpha q)X_{23}\otimes
X_{12}k^{-2}_2,\cr
{\rm and}\hskip0.2cm&\cr
&{\Delta X_{31}{=}\varepsilon \otimes X_{31}{+}X_{31}\otimes
(k_1k_2)^{-2}{+}(1{-}\beta q^{-1})X_{21}\otimes X_{32}k^{-2}_1{+}
(q^{-1}{-}\beta)X_{32}\otimes X_{21}k^{-2}_2}~.\cr}
$$
That is, we have $\Delta X-\varepsilon\otimes X\in\X\otimes{\cal{A}}'$ for all
$X\in\X$. Therefore, by Lemma 1, the above
Ansatz gives a left-covariant differential calculus $(\Gamma,\d)$
over ${\cal{A}}=SL_q(3)$.\nz
{}From the formulas for $\Delta X_{ij}$ and $\Delta X_n$ we see that
the corresponding homomorphism $f$ of the algebra $SL_q(3)$
(as defined in
Theorem 2.1, 3., in [18]) decomposes
into a direct sum of an upper triangular part, a lower triangular part
and a diagonal part. Using the pairing between $\U_q(sl_3)$ and 
$SL_q(3)$ we obtain
the following commutation relations between matrix entries and 1-forms:
$$\eqalign{
&\omega_{12} u^i_j=q^{-\delta_{j1}+\delta_{j2}} u^i_j\omega_{12}+\delta_{j3}
(q-\alpha)u^i_2\omega_{13},\cr
&\omega_{23} u^i_j=q^{-\delta_{j2}+\delta_{j3}} u^i_j\omega_{23}+\delta_{j2}
(q^{-1}-\alpha)u^i_1\omega_{13},\cr
&\omega_{13} u^i_j=q^{-\delta_{j1}+\delta_{j3}} u^i_j\omega_{13},\cr
&\omega_{21} u^i_j=q^{-\delta_{j1}+\delta_{j2}} u^i_j\omega_{21}+\delta_{j2}
(q-\beta)u^i_3\omega_{31},\cr
&\omega_{32} u^i_j=q^{-\delta_{j2}+\delta_{j3}} u^i_j\omega_{32}+\delta_{j1}
(q^{-1}-\beta)u^i_2\omega_{31},\cr
&\omega_{31} u^i_j=q^{-\delta_{j1}+\delta_{j3}} u^i_j\omega_{31},\cr
{\rm and}~\hskip2cm & \omega_n u^i_j=q^{-2\delta_{nj}+2\delta_{n+1,j}}
u^i_j\omega_n~.\hskip6cm\cr}
$$
In particular, if $\alpha=\alpha (q)$ and $\beta=\beta(q)$ are functions
of $q$ such that their limits are equal to $1$ as $q\rightarrow 1$, 
then the
calculus $(\Gamma,\d)$ gives the ordinary differential calculus on
$SL(2)$
in the limit $q\rightarrow 1$.
\nz
{}From the general theory [18] we know that the right ideal $\R$ of $\ker
\varepsilon$ associated with a FODC plays a crucial role. For the calculus
$(\Gamma,\d)$ defined above the
right ideal $\R$ is generated by the following elements:
$$\eqalign{
&u^r_su^i_j~ {\rm for}~r\ne s, r\ne j, i\ne j, i\ne s; u^r_ru^i_j-u^i_j ~
{\rm for}~i\ne j;\cr
&u^i_ju^j_i~ {\rm for}~i\ne j;~ u^2_1 u^1_3; u^3_1 u^1_2; u^1_3 u^3_2;
u^2_3 u^3_1;\cr
&u^2_3u^1_2 - (q^{-1}-\alpha)u^1_3; u^2_1u^3_2 - (q-\beta)u^3_1;\cr
&u^1_1u^3_3- u^1_1 - u^3_3+1; u^2_2u^1_1+q^{-2}u^3_3-(q^{-2}+1)1;
u^3_3u^2_2-u^2_2-q^{-2}u^3_3+q^{-2}1;\cr
&u^1_1u^1_1 - (1+q^{-2})u^1_1+q^{-2}1; u^2_2u^2_2-(1+q^{-2})u^2_2+
(q^4-1)u^1_1-(q^4-1-q^{-2})1;\cr
&u^3_3u^3_3 - (1+q^{-2})u^3_3+q^21;U-\varepsilon (U)1~.\cr}
$$
To prove this, we first verify by direct computations that the functionals
$X_{ij}$ and $X_n$ annihilate all these elements. We omit the (boring)
details of these computations. Thus $\R$ is contained in the right ideal of
$\ker\varepsilon$ associated with $(\Gamma,\d)$. Let $\E$ be the eight
dimensional vector space spanned by the elements $u^i_j$ for $i\ne j,
u^1_1-1$ and $u^2_2-1$. {}From the above list of generators of
$\R$ we conclude that each quadratic term $u^r_su^i_j-\delta_{rs}
\delta_{ij}1$ belongs to $\E + \R$. (For two missing terms
$u^1_2u^2_3$ and $u^3_2u^2_1$ we get $u^1_2u^2_3-(q-\alpha)u^1_3\in\R$
and $u^3_2u^2_1-(q^{-1}-\beta)u^3_1\in \R.)$ Hence we have ${\rm codim}~
\R\le 8$.
Since $\dim~ (\Gamma,\d)=8,\R$ is the right ideal of $\ker\varepsilon$
associated with the calculus $(\Gamma,\d)$.\nz
Now we turn to the higher order calculus of $(\Gamma,\d)$.
Our first aim is to prove the following
\bn
{\bf Lemma 4.} {\it If $(\alpha,\beta)\ne
(q,q^{-1})$ and $(\alpha,\beta)\ne (q^{-1},q)$, then we have $\omega_{13}
\otimes\omega_{31}\in\S1(\R)$.}\bn
{\bf Proof.} Recall that the elements $u^1_2u^2_1, u^2_3u^3_2$ and $u^1_3u^3_1$ belong to
the right ideal $\R$. We
determine the corresponding symmetric elements. All functionals $X_{ij}
X_{rs}$ with $i<j,r<s$ or $i>j, r>s$ and all functionals $X_nX_{rs},
X_{rs}X_n,X_nX_m$ annihilate these elements. We have
$$\eqalign{
\S1 (u^1_2u^2_1)&=X_{12} (u^1_2u^2_2)X_{21} (u^2_2u^2_1) \omega_{12}
\otimes\omega_{21}+X_{13}(u^1_3u^2_2)X_{31}(u^3_2u^2_1)\omega_{13}\otimes
\omega_{31} \cr
&+ X_{13} (u^1_2u^2_3) X_{31} (u^2_1u^3_2) \omega_{13}
\otimes\omega_{31}+X_{21} (u^1_1u^2_1) X_{12} (u^1_2u^1_1)
\omega_{21}\otimes\omega_{12}\cr
&=q \omega_{12} \otimes \omega_{21} + (q^{-1}-\beta) \omega_{13}
\otimes\omega_{31} + (q-\alpha) \omega_{13} \otimes \omega_{31}
+q^{-1} \omega_{21}\otimes \omega_{12}\cr
&=q \omega_{12} \otimes \omega_{21} + q^{-1} \omega_{21} \otimes
\omega_{12} + (\lambda_+-\alpha-\beta)\omega_{13} \otimes \omega_{31}.\cr}
$$
Similarly, we obtain
$$\eqalign{
&\S1 (u^2_3u^3_2) =q \omega_{23} \otimes \omega_{32} +q^{-1}
\omega_{32} \otimes \omega \omega_{23} + (\lambda_+-\alpha-\beta)
\omega_{31} \otimes \omega_{13}\cr
{\rm and}\hskip0.5cm &\S1 (u^1_3u^3_1)=q\omega_{13}\otimes
\omega_{31} +q^{-1}\omega_{31} \otimes \omega_{13}.\cr}
$$
Using two of these formulas and the facts that $\S1 (xy)=P_{inv} (\S1
(x)y),x\in\R$, by Lemma 1, and $P_{inv} (y \eta)=\varepsilon (y) P_{inv}
(\eta)$ for
$\eta\in\Gamma\otimes_{{\cal {A}}} \Gamma$ and $y\in{\cal A}$, we compute
$$\eqalignno{
&\S1 (u^1_2u^2_1u^3_3)=P_{inv} (\S1 (u^1_2u^2_1)u^3_3)\cr
&= P_{inv} (q \omega_{12} \otimes \omega_{21} u^3_3 + q^{-1}
\omega_{21} \otimes \omega_{12} u^3_3 + (\lambda_+-\alpha-\beta)\omega_{13}
\otimes \omega_{31} u^3_3 )\cr
&=P_{inv}(qu^3_3 \omega_{12}\otimes \omega_{21}+q(q-\alpha)u^3_2
\omega_{12}\otimes \omega_{21}+q^{-1} u^3_3 \omega_{21}
\otimes \omega_{12}\cr
&\hskip0.5cm+q^{-1}(q-\alpha) u^3_2 \omega_{21} \otimes \omega_{13}+q^{-1}
(q-\alpha)(q-\beta) u^3_3 \omega_{31}\otimes \omega_{13}+\cr
&\hskip0.5cm(\lambda_+-\alpha-\beta)q^2 u^3_3 \omega_{13}\otimes \omega_{31})\cr
&=q\omega_{12}\otimes \omega_{21} +q^{-1} \omega_{21}
\otimes \omega_{12} +q^{-1}(q-\alpha)(q-\beta)\omega_{31}
\otimes \omega_{13}+\cr
&\hskip0.5cmq^2(\lambda_+-\alpha-\beta)\omega_{13}\otimes \omega_{31}\cr
&=\S1 (u^1_2u^2_1)+(q-\alpha)(q-\beta)\S1 (u^1_3u^3_1)+(1-q\alpha)
(\beta-q^{-1})\omega_{13}\otimes \omega_{31},\cr}
$$
so that $\omega_{13}\otimes \omega_{31} \in\S1(\R)$
provided that $\alpha\ne q^{-1}$ and $\beta\ne q^{-1}$.\nz
Similarly, we get
$$\S1 (u^2_3u^3_2u^1_1)=\S1 (u^2_3u^3_2)+ (q^{-1}-\alpha)(q^{-1}-
\beta)\S1 (u^1_3u^3_1)+(\alpha-q)(\beta-q)
\omega_{13}\otimes\omega_{31},
$$
hence $\omega_{13}\otimes\omega_{31}\in\S1(\R)$ when
$\alpha \ne q$ and $\beta\ne q$.~\qed\sn\vfill\eject\noindent
Therefore, if $(\alpha,\beta)\ne (q^{-1},q)$ and $(\alpha,\beta)\ne
(q,q^{-1})$, then the 2-form $\omega_{13}\wedge \omega_{31}$ vanishes
in $\Gamma^\wedge_2$ and hence the space of 2-forms does {\it not} give the
corresponding space for the
ordinary differential calculus on $SL(3)$ when $q\rightarrow 1$.
The two remaining distinguished cases $(\alpha,\beta)=(q^{-1},q)$ and
$(\alpha,\beta)=(q,q^{-1})$ will be treated in the following section .
For this let $(\Gamma_1, \d)$ and $(\Gamma_2,\d)$ denote the first order
differential calculus $(\Gamma,\d)$ on $SL_q(3)$ with $(\alpha,\beta)=
(q^{-1},q)$ and $(\alpha,\beta)=(q,q^{-1})$, respectively.
\bn
{\bf 4. The differential calculi ${\bf(\Gamma_1,\d)}$ and
${\bf(\Gamma_2,\d)}$ on ${\bf SL_q(3)}$}
\mn
Let $\R_r$ be the right ideal of $\ker\varepsilon$ and let $\X_r$ be the
linear span of linear functionals $X_{ij},i\ne j,
i,j=1,2,3$, and $X_n,n=1,2$, for the calculus $(\Gamma_r,\d), r=1,2$.
{}From the definitions of functionals $X_{i,j},X_n$ and the commutation
rules in the algebra $\U_q(sl_3)$ we obtain the following commutation
relations for the generators of the quantum Lie algebra $\X_r$.
$$\eqalignno{
\X_1~{\rm and}~ \X_2:~
&X_{12}X_{32}{-}q^{-1}X_{32}X_{12} {=}0,X_{23}X_{21}{-}q^{-1}
X_{21}X_{23}{=}0~,\cr
&X_{12}X_{21}{-}q^{2}X_{21}X_{12} {=}X_1,X_{13}X_{31}{-}q^{2}
X_{31}X_{13}{+}q^{-1}\lambda X_1X_2{=}X_1{+}X_2~,\cr
&X_{23}X_{32}{-}q^{2} X_{32}X_{23} {=}X_2,X_{1}X_{2}{-}
X_{2}X_{1}{=}0~,\cr
&X_{1}X_{12}{-}q^{-4}X_{12}X_{1} {=}(1{+}q^{-2})X_{12},X_{1}
X_{21}{-}q^{4}X_{21}X_{1}{=}{-}(q^2{+}q^4)X_{21}~,\cr
&X_{2}X_{23}{-}q^{-4} X_{23}X_{2} {=}(1{+}q^{-2})X_{23},X_{2}
X_{32}{-}q^{4}X_{32}X_{2}{=}{-}(q^2{+}q^4)~X_{32},\cr
&X_{1}X_{23}{-}q^2 X_{23}X_{1} {=}{-}q^2X_{23},X_{1}X_{32}{-}
q^{-2}X_{32}X_{1}{=}X_{32}~,\cr
&X_{1}X_{13}{-}q^{-2} X_{13}X_{1} {=}X_{13},X_{1}X_{31}{-}
q^2X_{31}X_{1}{=}{-}q^2X_{31}~,\cr
&X_{2}X_{13}{-}q^{-2} X_{13}X_{2} {=}X_{13},X_{2}X_{31}{-}
q^2X_{31}X_{2}{=}{-}q^2X_{31}~,\cr
&X_{2}X_{12}{-}q^2 X_{12}X_{2} {=}{-}q^2X_{12},X_{2}X_{21}{-}
q^{-2}X_{21}X_{2}{=}X_{21}~,\cr
&{\ }\cr
\X_1:\hskip1.3cm
&X_{13}X_{23}{-}q X_{23}X_{13} {=}0,X_{32}X_{31}{-}q^{-1}
X_{31}X_{32}{=}0~,\cr
&X_{12}X_{13}{-}q X_{13}X_{12} {=}0,X_{31}X_{21}{-}q^{-1}
X_{21}X_{31}{=}0~,\cr
&X_{12}X_{23}{-}q^{-1} X_{23}X_{12} {=}X_{13},X_{32}X_{21}{-}
qX_{21}X_{32}{=}X_{31}~,\cr
&X_{12}X_{31}{-}q X_{31}X_{12} {=}{-}qX_{32},X_{13}X_{21}{-}q
X_{21}X_{13}{-}\lambda X_{23}X_1{=}{-}qX_{23}~,\cr
&X_{13}X_{32}{-}q X_{32}X_{13} {=}X_{12},X_{23}X_{31}{-}qX_{31}
X_{23}{+}q^{-1}\lambda X_{21}X_2{=}X_{21}~.\cr
{\ }\cr
\X_2:\hskip1.3cm
&X_{13}X_{23}{-}q^{-1} X_{23}X_{13} {=}0,X_{32}X_{31}{-}
qX_{31}X_{32}{=}0~,\cr
&X_{12}X_{13}{-}q^{-1} X_{13}X_{12} {=}0,X_{31}X_{21}{-}q
X_{21}X_{31}{=}0~,\cr
&X_{12}X_{23}{-}q X_{23}X_{12} {=}X_{13},X_{32}X_{21}{-}q^{-1}
X_{21}X_{32}{=}X_{31}~,\cr
&X_{12}X_{31}{-}q X_{31}X_{12} {-}\lambda X_{1}X_{32}{=}{-}qX_{32},
X_{13}X_{21}{-}qX_{21}X_{13}{=}{-}qX_{23}~,\cr
&X_{13}X_{32}{-}q X_{32}X_{13}{+}q^{-1}\lambda X_{2}X_{12}{=}X_{12},
X_{23}X_{31}{-}qX_{31}X_{23}{=}X_{21}.\cr}
$$
Next we shall desribe the corresponding higher order differential calculi.
For this we need to know the vector space $\S1(\R_r)$. Let I denote the
ordered index set $\{1,2,21,31,32,12,13,23\}.$
The right ideal $\R_r,r{=}1,2$, has 37 generators which have been listed in
the preceding section (recall that $\alpha=q^{-1},\beta=q$ for $\R_1$ and
$\alpha=q,\beta= q^{-1}$ for $\R_2$).
Let $\B_r$ be the linear span of these generators.
Using the formulas for the
comultiplications of $X_{ij}$ and $X_n$ and for the pairing between
$\U_q(sl_3)$ and $SL_q(3)$ one can compute the symmetric elements $\S1(x)$
for the generators of $\R_r$. We state only the result of these (long)
computations.
The following 36 elements of $\Gamma\otimes_{{\cal A}} \Gamma$ belong to
$\S1(\R_r)$ and form a basis of the vector space
$\S1({\cal B}_r)$.
$$\eqalign{
\S1 (\R_1)~{\rm and}~ \S1(\R_2):~
&\omega_{ij}\otimes \omega_{ij}~{\rm for}~ i\ne j, i,j{=}1,2,3;\omega_n\otimes \omega_n~
{\rm for}~n{=}1,2;\cr
&\omega_{12}\otimes \omega_{32}{+}q\omega_{32}\otimes\omega_{12},\omega_{23}
\otimes \omega_{21}{+}q\omega_{21}\otimes \omega_{23},\cr
&\omega_{12}\otimes \omega_{21}{+}q^{-2}\omega_{21}\otimes\omega_{12},\omega_{13}
\otimes \omega_{31}{+}q^{-2}\omega_{31}\otimes \omega_{13},\cr
&\omega_{23}\otimes \omega_{32}{+}q^{-2}\omega_{32}\otimes\omega_{23},\omega_{1}
\otimes \omega_{2}{+}\omega_2\otimes\omega_1-q^{-1}\lambda\omega_{13}\otimes
\omega_{31},\cr
&\omega_{1}\otimes \omega_{12}{+}q^4\omega_{12}\otimes\omega_{1},\omega_{1}
\otimes \omega_{21}{+}q^{-4}\omega_{21}\otimes \omega_{1},\cr
&\omega_{2}\otimes \omega_{23}{+}q^4\omega_{23}\otimes\omega_{2},\omega_{2}
\otimes \omega_{32}{+}q^{-4}\omega_{32}\otimes \omega_{2},\cr
&\omega_{1}\otimes \omega_{23}{+}q^{-2}\omega_{23}\otimes\omega_{1},\omega_{1}
\otimes \omega_{32}{+}q^2\omega_{32}\otimes \omega_{1},\cr
&\omega_{1}\otimes \omega_{13}{+}q^2\omega_{13}\otimes\omega_{1},\omega_{1}
\otimes \omega_{31}{+}q^{-2}\omega_{31}\otimes \omega_{1},\cr
&\omega_{2}\otimes \omega_{13}{+}q^2\omega_{13}\otimes\omega_{2},\omega_{2}
\otimes \omega_{31}{+}q^{-2}\omega_{31}\otimes \omega_{2},\cr
&\omega_{2}\otimes \omega_{12}{+}q^{-2}\omega_{12}\otimes\omega_{2},\omega_{2}
\otimes \omega_{21}{+}q^2\omega_{21}\otimes \omega_{2},\cr
\cr
\S1(\R_1):\hskip2cm
&\omega_{13}\otimes \omega_{23}{+}q^{-1}\omega_{23}\otimes\omega_{13},\omega_{32}
\otimes \omega_{31}{+}q\omega_{31}\otimes \omega_{32},\cr
&\omega_{12}\otimes \omega_{13}{+}q^{-1}\omega_{13}\otimes\omega_{12},\omega_{31}
\otimes \omega_{21}{+}q\omega_{21}\otimes \omega_{31},\cr
&\omega_{12}\otimes \omega_{23}{+}q\omega_{23}\otimes\omega_{12},\omega_{32}
\otimes \omega_{21}{+}q^{-1}\omega_{21}\otimes \omega_{32},\cr
&\omega_{12}\otimes \omega_{31}{+}q^{-1}\omega_{31}\otimes\omega_{12},\omega_{23}
\otimes \omega_{1}{+}q^2\omega_{1}\otimes \omega_{23}{+}\lambda\omega_{13}\otimes
\omega_{21},\cr
&\omega_{13}\otimes \omega_{32}{+}q^{-1}\omega_{32}\otimes\omega_{13},\omega_{2}
\otimes \omega_{21}{+}q^2\omega_{21}\otimes \omega_{2}{+}\lambda\omega_{31}\otimes
\omega_{23}.\cr
\S1(\R_2):\hskip2cm
&\omega_{13}\otimes \omega_{23}{+}q\omega_{23}\otimes\omega_{13},\omega_{32}
\otimes \omega_{31}{+}q^{-1}\omega_{31}\otimes \omega_{32},\cr
&\omega_{12}\otimes \omega_{13}{+}q\omega_{13}\otimes\omega_{12},\omega_{31}
\otimes \omega_{21}{+}q^{-1}\omega_{21}\otimes \omega_{31},\cr
&\omega_{12}\otimes \omega_{23}{+}q^{-1}\omega_{23}\otimes\omega_{12},\omega_{32}
\otimes \omega_{21}{+}q\omega_{21}\otimes \omega_{32},\cr
&\omega_{1}\otimes \omega_{32}{+}q^2\omega_{32}\otimes\omega_{1}{+}\lambda\omega_{12}
\otimes \omega_{31},\omega_{13}\otimes \omega_{12}{+}q^{-1}\omega_{21}
\otimes\omega_{21},\cr
&\omega_{12}\otimes \omega_{2}{+}q^2\omega_{2}\otimes\omega_{12}{+}\lambda\omega_{32}
\otimes \omega_{13},\omega_{23}\otimes \omega_{31}{+}q^{-1}\omega_{31}\otimes
\omega_{23},\cr}$$
The preceding formulas showed that both FODC $(\Gamma_1\d)$ and $(\Gamma_2,\d)$
are very close to the classical differential calculus on $SL(3)$. For
instance, except for two cases (that is, $\omega_{12} u^n_3,
\omega_{32} u^n_1~{\rm for}$
$r{=}1~{\rm and}~\omega_{23} u^n_2,\omega_{21} u^n_2~{\rm for}~
r{=}2)$, each 1-form $\omega_{{\goth{i}}} u^n_m$ is equal to $u^n_m
\omega_{{\goth{i}}}$ multiplied by
some power of $q$.
Except for three cases, the commutation relations of the quantum Lie algebra
$\X_r$ and the elements of the left module basis of $S_2$ contain only
two quadratic terms as in the classical case.
\bn
{\bf Lemma 5:} {\it For r{=}1,2, we have} ${\cal S}(\B_r){=}\S1 (\R_r)$.\mn
{\bf Proof.} By Lemma 2, (ii),
it suffices to show that $P_{inv}(\zeta u^i_j)\in\S1(\B_r)$ for $i,j{=}1,2,3$
and for all 36 basis elements $\zeta$ of listed above. This is obviously
fulfilled if $\zeta u^i_j$ is a scalar multiple of $u^i_j\zeta$, since
$P_{inv}(a\zeta){=}\varepsilon (a) P_{inv}(\zeta)$ by Lemma 2.2 in [18]. {}From
the commutation relations between matrix entries and 1-forms we conclude that
it remains to check the condition $P_{inv}(\zeta u^i_j)\in\S1 (\B_r)$ for all
basis elements $\zeta$ of $\S1(\B_r)$ which are sums of three terms or
contain the 1-forms $\omega_{12}$ or $\omega_{32}$ in case $r{=}1$ resp.
$\omega_{23}$ or $\omega_{21}$ in case $r{=}2$. These are 14 elements $\zeta$
in either case $r{=}1$ and $r{=}2$. We carry out this verification for the
two elements $\zeta_1:=\omega_{12}\otimes\omega_{32}{{+}}q\omega_{32}
\otimes\omega_{12}$ and $\zeta_2 :=\omega_1\otimes\omega_2{{+}}\omega_2\otimes
\omega_1-q^{-1}\lambda\omega_{13}\otimes\omega_{31}$ in case $r{=}1$.
{}From $\omega_{12}\otimes \omega_{32} u^i_3{=}q u^i_3 \omega_{12}
\otimes \omega_{32}{+}q\lambda u^i_2\omega_{13}\otimes \omega_{32}$ and
$\omega_{32}\otimes \omega_{12} u^i_3{=}qu^i_3\omega_{32}\otimes \omega_{12}
{+}q^{-1}\lambda u^i_2\omega_{32}\otimes \omega_{13}$ it follows that $P_{inv}
(\zeta_1u^3_3){=}q\zeta_1$ and $P_{inv}(\zeta_1u^2_3){=}q\lambda (
\omega_{13}\otimes \omega_{32}{+}q^{-1}\omega_{32}\otimes \omega_{13})\in
\S1 (B_r)$. Moreover, $\zeta_1u^i_j$ is a multiple of $u^i_j\zeta_1$ for all
elements $u^i_j$ other than $u^3_3$ and $u^2_3$. For $\zeta_2$ we obtain
$\zeta_2u^i_j{=}q^{-2\delta_{j1}{+}2\delta_{j3}}u^i_j \zeta_2$, so the
condition is also valid for $\zeta_2$. The other 26 cases are treated in a
similar way.~~\qed\sn
Thus, by Lemma 5, the 36 basis elements for the vector space $\S1(\B_r)$ of
the above list form a free left module basis for $S_2{=}{\cal A}\S1(\R_r){\cal A}$.
This gives a precise description of the ${\cal A}$-subbimodule $S_2$ of
$\Gamma\otimes_{{\cal A}}\Gamma$.\nz
{}From the definitions of functionals $X_{ij},X_n$ and the pairing between
$\U_q(sl_3)$ and $SL_q(3)$ it follows that for both calculi
$$
\omega (u^1_1)=\omega_1, \omega(u^2_2)=\omega_2-q^2\omega_1,\omega(u^3_3)=
-q^2\omega_2~{\rm and}~\omega (u^i_j)=\omega_{ij}~{\rm for}~i\ne j~.
$$
Recall that for any left-covariant differential calculus over ${\cal A}$ we have
$\d\omega (a){=}\break
-\omega (a_{(1)})\wedge\omega (a_{(2)}), a\in{\cal A}$. {}From the
preceding we obtain the following {\it Maurer-Cartan formulas} for our
differential calculi.
$$\eqalignno{
(\Gamma_1,\d)~{\rm and}~(\Gamma_2,\d):~
&\d \omega_1{=} {-}\omega_{12}\wedge \omega_{21}{-}\omega_{13} \wedge \omega_{31},
\d \omega_2 {=} {-}\omega_{13} \wedge \omega_{31} {-} \omega_{23} \wedge
\omega_{32},\cr
&\d \omega_{13}{=} {-}\omega_{1}\wedge \omega_{13}{-}\omega_{2} \wedge
\omega_{13}{-}\omega_{12}\wedge \omega_{23},\cr
&\d \omega_{23}{=} q^2 \omega_{1}\wedge \omega_{23}{-}(1{+}q^{-2})\omega_{2}
\wedge \omega_{23}{+}q \omega_{13} \wedge \omega_{21},\cr
&\d \omega_{21}{=} (q^2{+}q^4) \omega_{1}\wedge \omega_{21}{-}\omega_{2} \wedge
\omega_{21}{-}\omega_{23} \wedge \omega_{31}. \cr
\cr
(\Gamma_1,\d):\hskip1.8cm
&\d \omega_{12}{=} {-}(1{+}q^{-2})\omega_{1}\wedge \omega_{12}{+}q^2 \omega_{2}
\wedge \omega_{12}{-} \omega_{13} \wedge \omega_{32},\cr
&\d \omega_{31}{=} q^2 \omega_{1}\wedge \omega_{31}{+}q^2\omega_{2}
\wedge \omega_{31}{+}q^{-1} \omega_{21} \wedge \omega_{32},\cr
&\d \omega_{32}{=} {-}\omega_1 \wedge \omega_{32}{+}(q^2{+}q^4) \omega_{2}\wedge
\omega_{32}{+}q\omega_{12} \wedge \omega_{31}. \cr
\cr
(\Gamma_2,\d):\hskip1.8cm
&\d \omega_{12}{=} {-}(1{+}q^{-2})\omega_{1}\wedge \omega_{12}{+}q^2 \omega_{2}
\wedge \omega_{12}{-} q^2\omega_{13} \wedge \omega_{32},\cr
&\d \omega_{31}{=} q^2 \omega_{1}\wedge \omega_{31}{+}q^2\omega_{2}
\wedge \omega_{31}{+}q \omega_{21} \wedge \omega_{32},\cr
&\d \omega_{32}{=} {-}\omega_1 \wedge \omega_{32}{+}(q^2{+}q^4) \omega_{2}\wedge
\omega_{32}{+}q^{-1}\omega_{12} \wedge \omega_{31}. \cr}
$$
In this and the next paragraph we omit the subindex $r$ which refers to
one of the claculi $(\Gamma_1,\d)$ and $(\Gamma_2,\d)$. We define a $64\times 64$-matrix
$\sigma{=}(\sigma^{\i\j}_{\m\n})_{\i,\j,\n,\m\in I}$ as follows. Consider elements
$\omega_\n\otimes\omega_\m{+}\mu\omega_\m\otimes\omega_\n$ and $\omega_\i
\otimes\omega_\j+\gamma\omega_\j\otimes\omega_\i+\delta\omega_\n\otimes\omega_\m$
(written in the order of the index set $I$, i.e. ${\goth{n}}<{\goth{m}}$ and
${\goth{i}}<{\goth{j}}$) of our vector space basis
of $\S1(\R)$. We set $\sigma_{\n\m}^{\m\n}{=}\mu,\sigma_{\m\n}^{\n\m}=
\mu^{-1},\sigma_{\i\j}^{\j\i}=\gamma,\sigma_{\j\i}^{\i\j}=\gamma^{-1},
\sigma_{\i\j}^{\n\m}=\delta$ and $\sigma_{\j\i}^{\m\n}=-\delta\mu\gamma^{-1}$.
The number $\sigma^{\i\i}_{\i\i},\i\in I$, are set equal to $1$ and the
remaining matrix entries are set zero. Then the 36 basis elements of
$\S1 (\R)$ are $\omega_\i\otimes\omega_\j+{\sum\limits_{\n,\m}}
\sigma_{\i\j}^{\n\m}\omega_\n\otimes\omega_\m,\i,\j\in I,\i\le \j.$ {}From the above
formulas we see that the
commutation relations of the quantum Lie algebra generators can be
expressed as
$X_\i X_\j-{\sum\limits_{\n,\n}}\sigma_{\n\m}^{\i\j}X_\n X_\m={\sum\limits_\k}
C^\k_{\i\j}X_\k,
\i\ne \j$, with certain coefficients $C^\k_{\i\j}\in\bbbc$.
The linear transformation $\sigma$ of
the 64-dimensional vector space $\Gamma_{inv}\otimes\Gamma_{inv} $ has
eigenvalues 1 and $-1$ with multiplicities 36 and 28, respectively.
Obvoiusly, $\sigma^2$ is the identity. However, $\sigma$ does
{\it not} satisfy the
braid relation $\sigma_{12}\sigma_{23}\sigma_{12}=\sigma_{23}
\sigma_{12}\sigma_{23}$ on $\Gamma_{inv}\otimes\Gamma_{inv}\otimes
\Gamma_{inv}$. (For instance, we have $\sigma_{12}\sigma_{23}\sigma_{12}
(\omega_1\otimes\omega_{23}\otimes\omega_{13})=q\omega_{13}\otimes
\omega_{23}\otimes\omega_1-q^{-2}\lambda\omega_{13}\otimes\omega_{21}\otimes
\omega_{13}$ and $\sigma_{23}\sigma_{12}\sigma_{23}(
\omega_{1}\otimes\omega_{23}\otimes\omega_{13})=q\omega_{13}\otimes
\omega_{23}\otimes\omega_1-\lambda\omega_{13}\otimes\omega_{21}
\otimes\omega_{13}$
for the calculus $(\Gamma_1,\d)$.)\nz
Let $J(\R)=\oplus_n
J_n(\R)$ be the two-sided ideal of the tensor algebra $\Gamma_{inv}
^\otimes$ of the vector space $\Gamma_{inv}$ which is generated by the set
$\S1(\R)$. Clearly, $(\Gamma^\wedge)_{inv}$ is (isomorphic to) the
quotient algebra $\Gamma_{inv}^\otimes/J(\R)=\oplus_n\Gamma_{inv}^
{\otimes n}/J_n(\R)$. The above basis elements of the vector space
$\S1(\R)$ form a Gr\"obner basis (see. e.g. [6] for this concept) of the
ideal $J(\R)$ with respect to the above ordering of the index set $I$.
We have checked this by using
the computer algebra system FELIX [1]. (One may also verify this assertion
by performing explicit calculations.) Therefore, we get a vector space basis
for
the quotient space $\Gamma^\oplus_{inv}/J(\R)$ by taking all monomials
in the 1-forms $\omega_{{\goth i}}$, ${\goth i}\in I$, which are not multiples of the leading term
of one of these Gr\"obner basis elements. {}From the special form of these
elements we see that this set of monomials is the same as in the case
$q=1$. Hence the dimension of the vector space $\Gamma^{\otimes n}_{inv}/
J_n(\R)$ is equal to the corresponding dimension $\left( 8\atop n\right)$ in the
classical case. Moreover, it follows that the associated higher order
calculi of both calculi $(\Gamma_1,\d)$ and $(\Gamma_2,\d)$ give the
ordinary higher order calculus on $SL(3)$ in the limit $q\rightarrow
1$.
\bn
{\bf 5. Left-covariant differential calculi on $SL_q(N)$}\mn
Let $L^+=({^+l}^i_j)$ and $L^-=({^-l}^i_j)$ be the $N\times N$ matrices of
linear functionals ${^+l}^i_j$ and ${^-l}^i_j$ on ${\cal A} := SL_q(N)$ as defined
in [5]. Recall that $L^\pm$ is uniquely determined by the properties that
$L^\pm :{\cal A}\rightarrow M_N(\bbbc)$ is a unital algebra homomorphism and
${^\pm l}^i_j (u^n_m)=p^{\mp 1}(\hat{R}^{\pm 1})^{in}_{mj}$ for $i,j,n,m{=}1,{\dots},
N$, where $p$ is an $N$-th root of $q$ and $\hat{R}$ is the $R$-matrix
of the quantum group $SL_q(N)$.\nz
In this section we define $N^2-1$ dimensional left-covariant FODC
$(\Gamma_1,\d)$ and $(\Gamma_2,\d)$ over ${\cal A}=SL_q(N)$. They generalize the
two calculi over $SL_q(N)$ studied in the preceding section. For this we set
$$\eqalign{
&X_{ij}=\lambda^{-1}\kappa ({^-l}^j_i){^-l}^i_i~{\rm and}~X_{ji}=-\lambda^{-1}\kappa
({^+l}^i_j){^+l}^j_j~{\rm for}~i<j~{\rm and}~r=1,\cr
&X_{ij}:=-\lambda^{-1}{^+l}^j_j {^-l}^j_i~{\rm and}~X_{ji}=
\lambda^{-1}{^-l}^i_i{^+l}^i_j~{\rm for}~i<j~{\rm and}~r=2\cr
&X_n=q\lambda^{-1} (\varepsilon -({^-l}^n_n {^+l}^{n+1}_{n+1})^2)~{\rm for}~ n=1,
\dots,N-1~{\rm and}~ r=1,2.\cr}
$$
Let $\X_r$ denote the linear span of functionals $X_{ij},i\ne j,~i,j=1,
\dots, N$, and $X_n,n=1,\dots, N-1.$\nz
Computing the comultiplications of these generators, we obtain for $r=1$ and
$i<j$,\sn
\sn
$$
\Delta X_{ij}=\varepsilon \otimes X_{ij} +
{\sum\limits^{j}_{m=i+1}} X_{im}\otimes {^-l}^i_i\kappa ({^-l}^j_m)~
{\rm and}~\Delta
X_{ji}=\varepsilon\otimes X_{ji}+{\sum\limits^{j-1}_{m=i}} X_{jm}\otimes
{^+l}^j_j\kappa({^+l}^i_m)
$$
for $r=2$ and $i<j$,
$$
\Delta X_{ij}=\varepsilon \otimes X_{ij} +
{\sum\limits^{j-1}_{m=i}} X_{mj}\otimes {^+l}^j_j{^-l}^m_i~ {\rm and}~\Delta
X_{ji}=\varepsilon\otimes X_{ji}+{\sum\limits^j_{m=i+1}} X_{mi}\otimes
{^-l}^i_i {^+l}^m_j
$$
and for $r=1,2$,
$$
\Delta X_n=\varepsilon \otimes X_n+X_n
\otimes ({^-l}^n_n{^+l}^{n+1}_{n+1})^2~.
$$
In particular, these formulas show that $\Delta X-\varepsilon \otimes X\in
\X_r\otimes{\cal A}'$ for all $X\in\X_r$ and $r=1,2$. Therefore, by Lemma 1,
the
vector space $\X_r$ defines a left-covariant FODC over $SL_q(N)$. Furthermore,
we verify that $X_{ij}(u^r_s)=\delta_{ir}\delta_{js}$ for $i\ne j$ and
$X_n(u^r_s)=\delta_{rs}(\delta_{nr} -q^2\delta_{n+1,r})$ for $n=1,\dots,
N-1$. Thus all elements of the quantum Lie algebra $\X_r$ annihilate the quantum
trace $U=\sum q^{-2i} u^i_i$.
Since dim $\X_r=N^2-1$, the FODC $(\Gamma_r,\d)$ has dimension $N^2-1$.\nz
It is not difficult to check that in the classical limit $q\rightarrow 1$ both
FODC $(\Gamma_r,\d)$ give the ordinary FODC on $SL(N)$. (As in the preceding
sections, we define the limit $q\rightarrow 1$ of the calculus $(\Gamma_r,\d)$
by keeping the Maurer-Cartan basis $\omega (u^r_s), (r,s)\ne (N,N),$ of
$(\Gamma_r)_{inv}$ fixed.) Some computations show that for $|i{-}j|\ge 3$ and
$r=1,2$ the 2-form $\omega_{ij} \wedge\omega_{ji}$ vanishes in $(\Gamma_r)^\wedge_2$.
That is, if $N\ge 4$, both associated higher order calculi do {\it not} have
the classical higher order calculus on $SL(N)$ as their limits when $q
\rightarrow 1$. To overcome this disadvantage, we have taken up another
approach in [14].
\bn\bn
{\bf References}\mn
\litem{[1]} Apel, J. and Klaus, U., FELIX-an assistent for algebraists.
In: S.M. Watt (ed.), ISSAC' 91, pp.382-389, ACM Press, 1991
\sn
\litem{[2]} Connes, A., Non-commutative differential geometry, Publ. Math.
IHES {\bf 62}, 44-144 (1986)
\sn
\litem{[3]} Cuntz, J. and Quillen, D., Algebraic extensions and
nonsingularity, J. Amer. Math. Soc. {\bf 8}, 251-289 (1995)
\sn
\litem{[4]} Drinfeld, V.G., Quantum groups, Proceedings ICM 1986,
pp. 798-820, Amer. Math. Soc., Providence, 1987
\sn
\litem{[5]} Faddeev, L.K., Reshetikhin, N.Y. and Takhtajan, L.A.,
Quantization of Lie groups and Lie algebras, Algebra and Analysis {\bf 1},
178-206 (1987)
\sn
\litem{[6]} Geddes, K.O., Czapor, S.R. and Labahn, G.,
Algorithms for Computer Algebra, Kluwer Acad. Publ., Boston, 1992
\sn
\litem{[7]} Jimbo, M., A $q$-difference analogue of ${\cal U}(g)$
and the Yang-Baxter equation, Lett. Math. Phys. {\bf 22}, 177-186 (1991)
\sn
\litem{[8]} Lyubashenko, V. and Sudbery, A., Quantum supergroups of $GL(n,m)$
type: differential forms, Koszul complexes and Berezinians, Preprint,
York, 1993
\sn
\litem{[9]} Maltsiniotis, G., Le langage des espaces et de groups quantiques,
Commun. Math. Phys. {\bf 151}, 275-302 (1993)
\sn
\litem{[10]} Montogomery, S. and Smith, S.P., Skew derivations and
$\U_q(sl(2))$, Israel J. Math. {\bf 72}, 158-166 (1990)
\sn
\litem{[11]} M\"uller-Hoissen, F. and Reuten, C., Bicovariant differential
calculus on $GL_{p,q}(2)$ and quantum subgroups, J. Phys. A {\bf 26},
2955-2975 (1993)
\sn
\litem{[12]} Schm\"udgen, K. and Sch\"uler, A., Classification of bicovariant
differential calculi on quantum groups of type A,B,C and D, Commun.
Math. Phys. {\bf 167}, 635-670 (1995)
\sn
\litem{[13]} Schm\"udgen, K. and Sch\"uler, A., Classification of bicovariant
differential calculi on quantum groups, Commun. Math. Phys. {\bf 170},
315-336 (1995)
\sn
\litem{[14]} Schm\"udgen, K. and Sch\"uler, A., Left-covariant differential
calculi on quantum groups, to appear.
\sn
\litem{[15]} Stachura, P., Bicovariant differential calculus on
$S_\mu U(2)$, Lett. Math. Phys. {\bf 25}, 175-188 (1992)
\sn
\litem{[16]} Sudbery, A., Non-commuting coordinates and differential
operators, in: Quantum Groups, T. Curtright, D. Fairlie and C. Zachos (eds.),
pp. 33-52, World Scientific, Singapore, 1991
\sn
\litem{[17]} Woronowicz, S.L., Twisted $SU(2)$ group. An example of a
non-commutative differential calculus, Publ. RIMS Kyoto Univ. {\bf 23},
177-181 (1987)
\sn
\litem{[18]} Woronowicz, S.L.,  Differential calculus on compact matrix
pseudogroups (quantum groups), Commun. Math. Phys. {\bf 122}, 125-170 (1989)

\end